\input amstex.tex
\documentstyle{amsppt}

\topmatter
\title Twisted K theory invariants  \endtitle
\author Jouko Mickelsson \endauthor 
\affil Department of Mathematics, University of Helsinki, FIN-00014, Finland, and 
Mathematical Physics, Royal Institute of Technology, 
SE-106 91, Stockholm, Sweden \endaffil 
\date January 13, 2004 \enddate 
\endtopmatter 

\advance\vsize -2cm
\magnification=1200
\hfuzz=30pt

\define\tr{\text{tr}}
\define\gm{\bold{g}}
\document

\redefine\a{\alpha}
\redefine\b{\beta} 
\redefine\g{\gamma} 

\baselineskip= 18pt

\bf Abstract \rm  An invariant for twisted K theory classes on a
3-manifold is introduced. The invariant is then applied to the 
twisted equivariant classes arising from the supersymmetric Wess-Zumino-Witten 
model based on the group $SU(2).$  
It is shown that the classes defined by different highest 
weight representations of the loop group $LSU(2)$ are inequivalent. 
The results are compatible with Freed-Hopkins-Teleman identification 
of twisted equivariant K theory as the Verlinde algebra. 

\vskip 0.3in

\bf 0. Introduction \rm 

Twisted K theory classes arise in a natural way in two dimensional 
conformal field theory and can be described in terms of Verlinde
algebra, [FHT1-2], [AtSe]. In this paper I shall describe a rather elementary
method for a construction of numerical invariants for the twisted 
K theory classes in the case of an oriented, connected, simply
connected 3-manifold. The method is then applied to the case of the 
group manifold $SU(2)$ and it is shown that indeed the result matches 
the prediction in [FHT1-2]. The computations are based on the construction
of twisted K theory classes in terms of the supersymmetric Wess-Zumino-Witten 
model as described in [M]. The result, Theorem 2, shows that indeed
the equivariant twisted K theory classes constructed from different 
highest $SU(2)$ weights are different.

Ordinary complex K theory on a space $X$ can be defined as the abelian group 
(with respect to direct sums of Hilbert spaces) of homotopy classes 
of maps from $X$ to Fredholm operators in a complex Hilbert space 
$H.$ There is a grading mod 2 in complex K theory. The group $K^0(X)$
is defined by using the space of all Fredholm operators in $H$ whereas 
$K^1(X)$ is defined with the help of self-adjoint Fredholm operators
which have both positive and negative essential spectrum. 

To define twisted K theory one needs as an input a principal $PU(H)$ 
bundle $P$ over $X.$ Here $PU(H)$ is the projective unitary group 
$PU(H)= U(H)/S^1$ in the Hilbert space $H.$ These principal bundles
are classified by $H^3(X,\Bbb Z);$ an element $\omega\in H^3(X,\Bbb
Z)$ is called the Dixmier-Douady class of the bundle $P$ and it plays 
the role of the (first) Chern class for circle bundles. A bundle $P$ 
is called a \it gerbe \rm over $X.$ Usually a gerbe is equipped with 
additional structure, the gerbe connection which is a Deligne 
cohomology class on $X$ with top form $\omega.$ 

Given $P$ we can define an associated vector bundle 
$$Q = P\times_{PU(H)} \Cal F, \tag0.1$$ 
where $\Cal F$ denotes the space of (self-adjoint) Fredholm operators 
in $H$ and the action of $PU(H)$ on $\Cal F$ is defined by
conjugation, [BCMMS]. The twisted K theory $K^*(X, \omega)$ is then the set of
homotopy classes of sections of the bundle $Q.$ It is again an abelian 
group with respect to direct sums. 

As in the case of ordinary K theory, it is sometimes useful to have an 
alternative equivalent definition. In the case of $K^1(X)$ one can
replace self-adjoint Fredholm operators by unitary operators using the 
trick in [AS].  First one can contract to space of (unbounded) self-adjoint Fredholm
operators (with positive and negative essential spectrum) to bounded
self-adjoint operators with essential spectrum at the points $\pm 1.$
Then one can map these operators to unitaries by $F\mapsto
g=-\exp(i\pi F).$ The operator $g$ belongs to the group $U_1(H)$ of
unitary operators such that $g-1$ is a trace-class operator. The
advantage with this method is that we can explicitly produce the
generators $H^*(U_1(H), \Bbb Z)$ as differential forms 
$$ w_{2k+1} = \alpha_{2k+1} \tr\, (g^{-1} dg)^{2k+1} \tag0.2$$ 
where $\alpha_{2k+1}$ is a normalization coefficient. The $K^1$  theory  
classes on $X$ are then classified, modulo torsion, by the pull-backs
of classes $w_{2k+1}$ with respect to a mapping $X\to U_1(H).$ 

In the case of twisted K theory we can use the same trick simply by 
replacing in (0.1) the space $\Cal F$ by $U_1(H);$ this gives an 
alternative definition for $K^1(X,\omega).$ The case $K^0(X,\omega)$ 
has to be dealt differently. There is a different unitary group
$U_{res}(H)$ which has the same homotopy type as the space of all
Fredholm operators in $H.$ To define $U_{res}(H)$ one needs a
polarization
$H=H_{+}\oplus H_-$ to a pair of infinite-dimensional subspaces and a 
grading operator $\epsilon, $ such that $H_{\pm}$ has grade $\pm 1.$ 
The group $U_{res}(H)$ consists the of unitaries $g$ such that
$[\epsilon,g]$ is compact. In fact, instead of compactness one can as
well  require that $[\epsilon,g]$ belongs to some fixed Schatten ideal 
$L_p$ of bounded operators $A$ such that $|A|^p$ is trace-class, with 
$1\leq p < \infty, $  [P]. 

This work was partially supported by the Erwin Schr\"odinger Institute
for Mathematical Physics in Vienna.  
I want to thank Alan Carey for many useful discussions. 

\vskip 0.3in

\bf 1. Invariants for twisted K theory classes over a 3-manifold \rm 

\vskip 0.3in

Let $M$ be an oriented compact connected 3-manifold. Fix a triangulation of $M$
by a finite number of  closed sets $\Delta_{\alpha}\subset M,$ where each
$\Delta_{\alpha}$ is parametrized by a standard 3-simplex
(tetraed). We may assume without restriction that when the closed sets 
$\Delta_{\alpha}$ are slightly extended to open sets $U_{\alpha}$ then
$\{U_{\alpha}\}$ is a good cover of $M,$ i.e., all the multiple 
intersections of the open sets are contractible. 
A gerbe over $M$ is given
in terms of transition functions $\phi_{\alpha\beta}: U_{\alpha}\cap
U_{\beta} \to PU(H).$ Here $H$ is a complex (in general, infinite-dimensional) 
Hilbert space. 
Since the open sets are contractible, we may
lift these functions to maps $\phi_{\a\b}:U_{\a\b}=U_{\a}\cap U_{\b} \to U(H).$ The lifts  satisfy
$$ \phi_{\a\b}(x)\phi_{\b\g}(x)\phi_{\g\a}(x)= f_{\a\b\g}(x),\tag1.1$$ 
where $f_{\a\b\g} : U_{\a\b\g} \to S^1.$ 
Denote by $\tau$ the Dixmier-Douady class of the gerbe, given by the
above system of local functions satisfying 
$$f_{\a\b\g} f^{-1}_{\a\b\eta} f_{\a\g\eta} f^{-1}_{\b\g\eta} =1 \tag1.2$$
on quadruple overlaps.
For the logarithms of these functions we get
$$a_{\a\b\g\eta}= \log f_{\a\b\g}-\log f_{\a\b\eta} +\log
f_{\a\g\eta}-\log f_{\b\g\eta}= 2\pi i n\tag1.3$$  
for some integer $n.$ The sum of $a_{\a\b\g\eta}$'s evaluated at the
vertices $\Delta_{\a\b\g\eta}$ is then equal to $2\pi i$ times the integral of the
Dixmier-Douady class over the 3-manifold $M.$ This can be written as 
$2\pi ik,$ where $k$ is an integer depending only on the Dixmier-Douady
class. 

A twisted $K^1$ theory class with a gerbe as input is then given by a family of
functions $g_{\a}: U_{\a} \to U_1(H).$ Here $U_1(H)$ is the group 
of unitaries $g$ in $H$ such that $g-1$ is trace-class. On the
overlaps $U_{\a\b}$,
$$ g_{\a} = \phi_{\a\b} g_{\b} \phi_{\a\b}^{-1}.\tag1.3$$ 

We want to determine a homotopy invariant for this class $[\omega]\in 
K^1(M, \tau).$ 

First let us recall that  an untwisted K theory class is the homotopy 
class of a globally defined function $g: M\to U_1(H).$ A homotopy
invariant for this is
$$I_M(g)= \frac{1}{24\pi^2} \int_M \tr\, (g^{-1} dg)^3,\tag1.4$$
i.e., the Witten action. This is an integer depending on the homotopy 
class of $g.$ 

In the twisted case we could try to use the formula 
$$\sum_{\a} \int_{\Delta_{\a}} \tr\, (g_{\a}^{-1} dg_{\a})^3.$$ 
However, this fails to be homotopy invariant due to boundary terms 
in integration by parts.    Instead, we can add correction terms 
$$r_{\a\b}= \int_{\Delta_{\a\b}}\omega_{\a\b}\tag1.5$$
with 
$$ \omega_{\a\b}=\frac{1}{8\pi^2}  \tr\, (d\phi_{\a\b}\phi_{\a\b}^{-1}) 
[ dg_{\a} g_{\a}^{-1} + g_{\a}^{-1}dg_{\a} + 
g_{\a} d\phi_{\a\b} \phi_{\a\b}^{-1}g_{\a}^{-1}-d\phi_{\a\b} \phi_{\a\b}^{-1} ].     
\tag1.6$$ 
Note that the second and the third term in the brackets are not
trace-class operators but their difference is.  
These correction terms are chosen such that
$$ d\omega_{\a\b}= \frac{1}{24\pi^2}\tr\, [(g_{\b}^{-1}dg_{\b})^3 -
(g_{\a}^{-1} dg_{\a})^3].
\tag1.7$$
 
Suppose for a moment that  all $f_{\a\b\g}=1.$ Then 
$$\omega_{\a\b}+\omega_{\b\g}+\omega_{\g\a}=0 \tag1.8$$
on triple overlaps of open sets. Define
$$I_M'(g)= \sum_{\alpha} I_{\Delta_{\a}}  
             + \sum_{\a< \b} r_{\a\b}  \tag1.9$$
where we have chosen the finite index set to be $\{1,2,\dots, p\}$ so
that we have a natural  ordering $\a < \b.$ Then it is a direct
consequence of Stokes' theorem, the Cech - de Rham cocycle relations
(1.7), (1.8), and 
closedness of  the forms $\tr\, (g^{-1}dg)^3$ that $I'_M(g)$ is a homotopy
invariant. 

However, in the case of a nontrivial gerbe the functions
$f_{\a\b\g}\neq 1$ and the cocycle relation (1.8) does not hold. 
The correct relation is 
$$\omega_{\a\b}+\omega_{\b\g} +\omega_{\g\a} = d\omega_{\a\b\g},\tag1.10$$ 
where $\omega_{\a\b\g}$'s are 1-forms on triple overlaps. A solution of
(1.10) is given by 
$$ \omega_{\a\b\g}= \frac{1}{4\pi^2}  h^{-1}dh\, \log\,f_{\a\b\g} , \tag1.11$$ 
where $h: M \to S^1$ is the globally defined function $h=\text{det}
g_{\a}.$ The choice of the index $\a$ is unimportant, since 
$g_{\a} = \phi_{\a\b} g_{\b} \phi_{\a\b}^{-1}$ so that the determinant
is well defined. However, we have to make a choice of the logarithm 
$\log\, f_{\a\b\g} .$ Two different choices differ by the locally constant 
function $n \cdot 2\pi i$ and give two different solutions to (1.9). 
In any case, the cocycle property (1.2) shows that 
$$\omega_{\a\b\g}
-\omega_{\a\b\eta}+\omega_{\a\g\eta}-\omega_{\b\g\eta}
= a_{\a\b\g\eta} h^{-1}dh
\tag1.12$$ 
on quadruple overlaps.  
Now we make the additional assumption that the function $h:M \to S^1$ 
is contractible (which would be automatic if $M$ is simply connected).
Then we can write
$$a_{\a\b\g\eta} h^{-1}dh= d(a_{\a\b\g\eta} \log\, h)
\equiv d\omega_{\a\b\g\eta}\tag1.14$$
with some choice of logarithm of $h.$ Different choices of the
logarithm lead to expressions for $\omega_{\a\b\g\eta}$ which 
differ by  $(2\pi i)^2$ times an integer. 

From this we reduce, by Stokes' theorem: 

\proclaim{Theorem 1} Let the determinant function $h$ defined above be
contractible. Then the expression
$$I(g)= \sum_{\a}  \int_{\Delta_{\a}} \omega_{\a} -\sum_{\a <\b}
\int_{\Delta_{\a\b}} \omega_{\a\b} + \sum_{\a<\b < \g}
\int_{\Delta_{\a\b\g}} \omega_{\a\b\g}
+\sum_{\a<\b<\g<\eta} \omega_{\a\b\g\eta} 
$$ 
is a homotopy invariant; the last term is evaluated at the points
 $\Delta_{\a\b\g\eta}.$
\endproclaim 

\bf Remark 1 \rm  The quantity $I(g)$ is only well defined modulo $k\times$
an integer. This is because of the arbitrary choice of the branch of
the logarithm of $h.$ The difference between two choices gives a
contribution 
$$\delta=2\pi i\cdot  \frac{1}{4\pi^2}\sum  a_{\a\b\g\eta}.$$ 
   
The  sum of the numbers $a_{\a\b\g\eta}$ is  equal to $2\pi i k \times$ an integer, where $k$ is 
an integer depending only on the Dixmier-Douady class $\tau$ of the
gerbe. Thus  
$\delta$ is equal to $k \times$ an integer and $I(g)$ is well defined mod $k.$ 

\bf Remark 2 \rm In the case when $h$ is not contractible we can still
use it to define the winding number invariant for the $K$ theory class,
$$w(h) = \frac{1}{2\pi i} \int_{S^1} h^{-1} dh,$$ 
where $S^1\subset M$ represents any element of $\pi_1(M).$

\bf Example \rm Take $M= S^3=SU(2).$ Then $H^3(M,\Bbb Z)$ is one
dimensional, the Dixmier-Douady class $\tau$ is represented as
$k$ times the basic 3-form $\frac{1}{24\pi^2} \tr\, (g^{-1}dg)^3$ 
on $SU(2).$ The map $I$ takes values in $\Bbb Z/ k\Bbb Z.$

\vskip 0.3in

\bf 2. Calculations in the case $G=SU(2)$ \rm 

\vskip 0.3in

We study the twisted K theory class over the group $G=SU(2). $ The Lie 
algebra of $G$ is denoted by $\gm.$ 
Let $\Cal A$ denote the space of smooth $\gm$ valued vector potentials
(1-forms) on the unit circle $S^1.$ Let $LG$ be the group of smooth
loops in $G$ and let $\Omega G\subset LG$ be the group of based loops,
i.e., loops $f$ such that $f(1)$ is the neutral element in $G.$ Then
$\Cal A/\Omega G$ is the group $G$ of holonomies around the
circle. The right action on $\Cal A$ is defined by $A^f = f^{-1} A f 
+f^{-1}df.$ 
The twisted K theory classes are constructed using the family
of  hermitean operators $Q_A$ for $A\in\Cal A$
constructed in [M]. 

The operator $Q_A$ is a sum of a 'free' supercharge $Q$ and an
interaction term $\hat A.$ The Hilbert space $H$ is a tensor 
product of a 'fermionic' Fock space $H_f$ and a 'bosonic' Hilbert 
space $H_b.$ The space $H_b$ carries an irreducible representation 
of the loop algebra $L\gm$ of level $k$ where The highest 
weight representations of level $k$ are classified by the $SU(2)$ 
representation of dimension $2j_0 +1$ on the 'vacuum sector'. 
We denote the generators of the loop algebra 
by $T_n^a,$ where $n\in\Bbb Z$ is the Fourier index and $a=1,2,3$ 
labels a basis of $\gm.$ The commutation relations are 
$$[T_n^a,T_m^b]= \lambda_{abc} T_{n+m}^{c} + \frac{k}{4} \delta_{ab}
\delta_{n,-m},\tag2.1$$ 
where $a,b,c=1,2,3$ are the structure constant of $\gm;$ in this case
when $\gm$ is the Lie algebra of $SU(2)$ the nonzero structure
constants are completely antisymmetric and we use the normalization 
$\lambda_{123}=
\frac{1}{\sqrt{2}}.$ (This comes from a normalization of the basis
vectors
$T^a_0\in \gm$ with respect to the Killing form.) 
In addition, we have the hermiticity relations
$(T^a_n)^*= -T^a_{-n}.$ With this normalization of the basis, $k$ is a
nonnegative integer and $2j_0= 0,1,2\dots k.$ The case $k=0$ corresponds
to a trivial representation and we shall assume in the following that
$k$ is strictly positive. 

The Fock space $H_f$ carries an irreducible representations of the
canonical anticommutation relations (CAR), 
$$\psi_n^a \psi_m^b +\psi_m^b\psi_n^a =
2\delta_{ab}\delta_{n,-m},\tag2.2$$ 
and $(\psi^a_n)^* = \psi^a_{-n}.$ The representation is fixed by the 
requirement that there is an irreducible representation of the
Clifford algebra $\{\psi_0^a\}$ in a subspace $H_{f,vac}$ such that 
$\psi_n^a v=0$ for $n<0$ and $v\in H_{f,vac}.$ 

The central extension of the loop algebra at level $2$ is
represented in $H_f$ through the operators 
$$K^a_n = -\frac14 \sum_{b,c=1,2,3; m\in\Bbb Z} \lambda_{abc} \psi_{n-m}^b
\psi_{m}^c, \tag2.3$$ 
that is, 
$$[K^a_n, K^b_m]= \lambda_{abc} K^c_{n+m} +\frac12 n \delta_{ab}\delta_{n,-m}.  
\tag2.4$$ 

We set $S^a_n= T^a_n + K^a_n.$ This gives a representation of the loop
algebra at level $k+2$ in the tensor product $H=H_f\otimes H_b.$ 

Next we define 
$$Q= i \psi^a_n T^a_{-n} +  \frac{i}{3} \psi^a_n K^a_{-n}.\tag2.5$$ 
This operator satisfies $Q^2 =h,$ where $h$ is the hamiltonian of the 
supersymmetric Wess-Zumino-Witten model, 
$$h= - \sum_{a,n} : T^a_n T^a_{-n} : + 
\frac{k+2}{8} \sum_{a,n} :n \psi^a_n\psi^a_{-n}: +\frac18,\tag2.6$$
where the normal ordering $::$ means that the operators with negative
Fourier index are placed to the right of the operators with positive
index,  $:\psi^a_{-n} \psi^b_{n}:\,\, = -\psi^b_{n} \psi^a_{-n}$ if $n>0$ 
and $:AB: = AB$ otherwise. In the case of the bosonic currents $T^a_n$ 
the sign is $+$ on the right-hand-side of the equation. See [KT] for
details on the supersymmetric current algebra. 

Finally, $Q_A$ is defined as 
$$Q_A = Q + i\tilde k \psi^a_n A^a_{-n}\tag 2.7$$ 
where the $A^a_n$'s are the Fourier components of the $\gm$-valued 
function $A$ in the basis $T^a_n$ and $\tilde k = \frac{k+2}{4}.$ 
All the formulas above can be generalized in a straight-forward way
to arbitrary simple Lie algebras, with the modification that the
last term $1/8$ in (2.6) is replaced by dim$\bold g/24$ and the level 
$k$ is quantized as integer times twice the lenght squared of the
longest root with respect to the dual Killing form.  

The basic property of the family of self-adjoint Fredholm operators
$Q_A$ is that it is equivariant with respect to the action of the 
central extension of the loop group $LG.$ Any element $f\in LG$ is 
represented by a unitary operator $S(f)$ in $H$ but the phase of $S(f)$ 
is not uniquely determined. The equivariantness property is
$$S(f^{-1}) Q_A S(f) = Q_{A^f}\tag2.8$$ 
with $A^f = f^{-1} A f +f^{-1} df.$    The infinitesimal version of
this is 
$$[S^a_n, Q_A] = i{\tilde k} ( n\psi^a_n + \sum_{b,c; m} \lambda_{abc}\psi^b_m A^c_{n-m})  
\tag2.9$$ 
which can be checked directly from (2.1), (2.2), and (2.4). 

The group $LG$ can be viewed as a subgroup of the group $PU(H)$
through the projective repsentation $S.$ The space $\Cal A$ of smooth 
vector potentials on the circle is the total space for a principal 
bundle with fiber $\Omega G \subset LG.$ Since now $\Omega G \subset PU(H),$ 
$\Cal A$ may be viewed as a reduction of a $PU(H)$ principal bundle 
over $G.$ The $\Omega G$ action by conjugation on the Fredholm
operators in $H$ defines an associated vector bundle $\Cal Q$ over 
$G$ and the family of operators $Q_A$ defines a section of this vector 
bundle. Thus $\{Q_A\}$ is a twisted $K^1$ theory class over $G$ where
the twist is determined by the level $k+2$ projective
representation of $LG.$

Using the method in [AS] we replace the family of unbounded hermitean
operators by a family of bounded operators $F_A= Q_A/(|Q_A| +
e^{-Q_A^2})$ which represent the same K theoretic class. The
perturbation in the denominators is introduced to avoid singularities 
with zero modes of $Q_A.$ The operator $F_A$ differs from the sign 
operator $Q_A/|Q_A|$ by a trace-class perturbation. For this reason 
the unitary operators $g_A= -e^{i\pi F_A}$ differ from the unit by a 
trace-class operator. 

We shall now study the twisted K theory class represented by the
family $g_A$ of unitary operators. Note that this family is still 
gauge equivariant, 
$$ S(f)^{-1} g_A S(f) = g_{A^f}\tag2.10$$ 
where $f\in LG.$ 

Since $S^3= \Cal A/\Omega G,$ we write the K theory class as a
function from the three dimensional unit disk $D^3$ to unitaries 
of the form $1+$ trace-class operators such that on the boundary 
$S^2$ the operators are gauge conjugate. Concretely, this is 
achieved as follows. For each point $\bold n\in S^2$ we define
a constant  $SU(2)$ vector potential $A(\bold n)=\frac{1}{2i} \bold n \cdot \bold
\sigma.$ 
Pauli matrices satisfy
$\sigma_1\sigma_2=i\sigma_3$ (and cyclic permutations) and
$\sigma_j^2=1.$  
The holonomy around the circle $S^1$ is equal to $-1$ for each
of the potentials $A(\bold n),$  thus they belong to the same $\Omega G$
orbit in $\Cal A.$ Next we define a disk $D^3$ of potentials $A(t,\bold n)=
t A(\bold n)$ where $0\leq t\leq 1$ is the radial variable in the disk $D^3$
and $\bold n$ are the angular coordinates. This disk projects to a closed
sphere $S^3$ in $G=\Cal A/\Omega G.$ For each $A\in D^3$ we have the
corresponding supercharge
$$Q_A= Q +  \frac{k+2}{4} t\cdot \sqrt{2}\psi^a_0  n^a\tag2.11$$
where the factor $\sqrt{2}$ comes from the normalization of the basis
$T^a$ of $\gm$ relative to the Pauli matrix basis.

Now we have a family of unitaries $g(t,\bold n)= g_{A(t,\bold n)}$ 
which are gauge conjugate on the boundary through the projective
unitary representation of $LG$ of level $k+2.$ This means that 
the homotopy class of the functions $g(t,\bold n)$ gives an element 
in ${K^1}_G(S^3, k+2).$ In the language of section 1, we may replace 
the triangulation $\{\Delta_{\a}\}$ by two sets: the disk $D^3$ as 
the southern hemisphere of $S^3$ and a second disk ${D'}^3$ as the
northern hemisphere. On the southern hemisphere we have the unitary 
matrix valued function $g(t,\bold n)$ whereas on the northern
hemisphere we have a constant function $g_0= -\exp(\pi i F_0).$ On the
equator parametrized by $\bold n\in S^2$ they are all gauge conjugate. 

The $G$ equivariantness follows from the fact
that the family $Q_A$ is gauge equivariant with respect to the full 
group $LG$ of gauge transformations and not only with respect to the 
based gauge transformations $\Omega G.$ 

We want to compute the quantum invariant for the class $[g]$ by 
evaluating the Witten functional 
$$I(g)= \frac{1}{24\pi^2} \int_{D^3} \tr (g^{-1} dg)^3.\tag2.12$$ 
Note that in the present setting the correction terms are absent since 
$g_0$ is constant (which can be deformed to the unit matrix since the 
group of unitaries $U_1(H)$ is connected).  
A direct computation of the integral of the trace in an
infinite-dimensional Hilbert space $H$ is difficult. Instead, we shall 
apply first various homotopy deformations to $g$ to bring the trace 
into more manageable form. 

\bf First deformation. \rm We need first a Lemma: 

\proclaim{ Lemma 1} The spectral projections $P_{\Lambda}$ of $|Q|$ 
commute with $Q_A$ when $A=\frac{1}{2i} t \bold n\cdot{\bold \sigma}.$
\endproclaim

\demo{Proof} Now $Q_A$ is given by (2.11). Using the canonical
anticommutation relations for $\psi_n^a$'s we observe that 
$$[Q_A-Q, Q]_+ =  -2\tilde k S^a_0 A_0^a.$$ 
On the other hand, $[S_0^a, Q]=0$ so that 
$$[Q_A, Q^2]= [Q_A -Q, Q^2]= 2\tilde k (- S^a_0 A_0^a Q + 
Q S^a_0 A_0^a) =0$$ 
from which follows $[Q_A, |Q|]=0$ and thus also $[Q_A, P_{\Lambda}]=0$ 
where $P_{\Lambda}$ is the spectral projection $|Q| \leq \Lambda.$  
\enddemo 

The Lemma implies that the spectral subspaces
$H_{\Lambda}=P_{\Lambda}H$ and $H_{\Lambda}^{\perp}$ are invariant
under $Q_A, F_A,$ and $g_A.$ 

Since $(Q_A -Q)^2 = 2 t^2 \tilde k^2$ we see that the restriction of 
$Q_A$ to the subspace $H_{\Lambda}^{\perp}$ is invertible if we choose 
$\Lambda > \sqrt{2} \tilde k.$ 

Let us deform the denominator $|Q_A|+ e^{{-Q_A}^2}$ in $F_A.$ Define
$D(s) = |Q(t,\bold n)| + s t(1-t) e^{-Q^2} + (1-s)e^{-Q(t,\bold
n)^2}$ for $0\leq s\leq 1.$ For any fixed $s$ these operators are 
gauge conjugate at the boundary $t=1$ because at $t=1$ we have 
$D(s)= |Q(1,\bold n)|+(1-s)e^{-Q(t,\bold n)^2}$. 
At $s=0$ this is the original family of 
denominators whereas for $s=1$ we get $D(1)= |Q(t,\bold n)| + 
t(1-t) e^{-Q^2}.$ 
This is our first deformation: We replace the original $g(t,\bold n)$ 
by the homotopy equivalent family 
$$g(t,\bold n)= -e^{i\pi F(t,\bold n)} \text{ with } F(t,\bold n)= 
\frac{Q(t,\bold n)}{|Q(t,\bold n)| + t(1-t) e^{-Q^2}}.$$

\bf Second deformation \rm By a similar $s$ dependent family as above 
we can replace the denominator $D(n)$ by $|Q(t,\bold n)|
+t(1-t)P_{\Lambda}$ for any $\Lambda > \tilde k.$ This is because
$Q(t,\bold n)$ is invertible in the complement of $H_{\Lambda}$ and 
is invertible in the whole space $H$ for $t=0,1.$ (For $t=0$ this is
clear since $Q^2\geq 1/8$ and for $t=1$ one observes that the spectrum
of $Q(1,\bold n)^2= Q^2+ \frac18 (k+2)^2 + i\sqrt{2}(k+2)\bold
n\cdot\bold S_0$ is of the form $\frac18 [1 + (k+2)p]$ where $p$ 
is an integer.)  
For the intermediate 
values $0<t<1$ both $t(1-t)e^{-Q^2}$ and $t(1-t)P_{\Lambda}$ are 
strictly positive in $H_{\Lambda}.$ So after the second deformation
$$g(t,\bold n)= -e^{i\pi F(t,\bold n)} \text{ with } F(t,\bold n)= 
\frac{Q(t,\bold n)}{ |Q(t,\bold n)| + t(1-t) P_{\Lambda}}.\tag2.13$$ 
In particular, since the eigenvalues of $Q/|Q|$ are $\pm 1,$ the
restriction of $g$ to $H_{\Lambda}^{\perp}$ is equal to the unit 
operator. Thus 
 
$$ I(g)= \frac{1}{24\pi^2} \int_{D^3} \tr_{H_{\Lambda}}
(g^{-1}dg)^3.\tag2.14$$
We need now to compute the trace of $(g^{-1}dg)^3$ only in the
finite-dimensional subspace $H_{\Lambda}.$ 
We use the formula 
$$\tr (g^{-1}dg)^3= d \tr\, dX \eta(ad_X)dX,\tag2.15$$ 
where $X=\log(g)$ and $\eta(x)= \frac{\sinh(x)-x}{x^2}.$ 
By Stokes' theorem the integral defining $I(g)$ is then equal 
to the integral of the 2-form $\tr\, dX \eta(ad_X) dX$ over 
$S^2= \partial D^3.$ But on the boundary $t=1$ we have $F(1,\bold n)= 
Q(1,\bold n)/|Q(1,\bold n)|.$ This simplifies $\eta(ad_X)$ so that 
the 2-form becomes
$$ \frac{i}{16\pi} \tr\, F dF dF  \text{ for } F=F(1,\bold n).$$ 
Summarizing  we obtain
$$I(g)=  \frac{i}{16\pi}  \int_{S^2} \tr_{H_{\Lambda}} F dF dF.\tag2.16$$

\bf Third deformation \rm We use the fact that the parameter $\Lambda$ 
is free except for the constraint $\Lambda > \sqrt{2}\tilde k.$ Since the
spectrum of $Q$ is discrete (the eigenvalues of $Q^2$ are quantized in units $(k+2)/2$), 
we can choose $\Lambda -\sqrt{2}\tilde k$ so 
small that the eigenvalues of $|Q|^2$ which are smaller or equal to
$\Lambda^2$ are also strictly smaller than $2 \tilde k^2.$ With this
choice $Q_A$ becomes invertible in $H_{\Lambda}.$ Furthermore,
also $Q_s(\bold n)=sQ +\sqrt{2} \tilde k \bold n\cdot \bold \psi_0$ is invertible in 
$H_{\Lambda}$ for all $0\leq s\leq 1.$ We use the homotopy $Q_s$ to 
replace $F= Q(1,\bold n)/|Q(1,\bold n)|$ in the integral $I(g)$ by 
the operator $F=\bold n\cdot \bold \psi_0.$  Now
$$\tr_{H_{\Lambda}} \, F dF dF =  {\bold n}\cdot
d\bold n\times d\bold n \,\tr_{H_{\Lambda}}\,\psi_0^1 \psi_0^2\psi_0^3.\tag2.17$$ 
The integral of $\bold n\cdot d\bold n\times d\bold n$ over $S^2$ 
is equal to twice the volume of $S^2$ and so
$$I(g)= - \frac{i}{2}\tr_{H_{\Lambda}} \, \Gamma  \text{ with } 
\Gamma= \psi_0^1\psi_0^2\psi_0^3.\tag2.18$$       
The trace  is essentially the Witten index. The operator $\Gamma$
almost anticommutes with the supercharge $Q.$ Define 
$$Q_+ = i\sum_{n\neq 0}\psi^a_n T^a_{-n} + \frac{i}{12} \sum_{n,m,n+m\neq 0}
\lambda_{abc} \psi^a_n\psi^b_m\psi^c_{-n-m}.\tag2.19$$
We have $Q_+ \Gamma = -\Gamma Q_+$ since $\psi^a_0$ anticommutes with 
$\psi^b_n$ for every $n\neq 0.$   

\proclaim{Lemma 2} The operator $Q_+$ commutes with the spectral
projections $P_{\Lambda}.$ \endproclaim
\demo{Proof} Write $Q= Q_0 + Q_+.$ Then $Q_0$ commutes with $\Gamma.$ 
Now 
$$Q^2 = Q_0^2 + Q_+^2 + [Q_0, Q_+]_+ =h$$ 
is even with respect to $\Gamma.$ The first two terms on the right are 
even, so the third term which is odd has to vanish and so $h=Q_0^2
+Q_+^2.$ This implies $[h, Q_+]= [Q_0^2, Q_+]= Q_0[Q_0,Q_+]_+ -
[Q_0,Q_+]_+ Q_0 = 0$ and so the spectral projections of $h$ commute 
with $Q_+.$ Since $|Q|^2=h,$ the same is true for the spectral
projections
$P_{\Lambda}$ of $|Q|.$ \enddemo 

\proclaim{Lemma 3} ${Q_+}^2 = Q^2 - \sum_a ({T^a}_0 + {{K'}^a}_0)^2-\frac{N}{24}$ where
${{K'}^a}_0 = -\frac14 \sum_{n\neq 0; b,c}\lambda_{abc} {\psi^b}_n {\psi^c}_{-n}.$
The ${{K'}^a}_0$'s satisfy the same commutation relations as the $T^a$'s.
\endproclaim \demo{Proof} By a direct computation. \enddemo

\proclaim{Lemma 4} Let $G=SU(2).$ Then the kernel of $Q_+$ in $H_{\Lambda}$
is equal to the vacuum sector
$H_0 \subset H$ consisting of eigenvectors of $h$ associated to the 
minimal eigenvalue $1/8.$  \endproclaim

\demo{Proof} Clearly $H_0\subset ker\, Q_+$ since $\psi^a_n v = T^a_n
v =0$ for any $v\in H_0$ for $n<0.$ We have to show that $|Q_+|$ is strictly 
positive in the orthogonal complement $H_0^{\perp}$ in $H_{\Lambda}.$ 

Let $d$ be the derivation in the affine Lie algebra based on $SU(2).$ 
By definition, $[d, T^a_n] = nT^a_n$ and $[d,\psi^a_n]= n\psi^a_n.$ 
From the weight inequalities for lowest  weight representations 
of affine Lie algebras, [K], Prop. 11.4, follows that in a $SU(2)$
subrepresentation with angular momentum $\ell,$ 
$$\ell_0(\ell_0+1)- d_0(k+2) \geq  \ell(\ell +1) - d(k+2),\tag2.20$$
where $\ell_0$ is the angular momentum of the lowest weight vector and
$d_0$ is the eigenvalue of $d$ for the lowest weight vector. 

We first apply the inequality to the bosonic representation in $H_b.$
The bosonic hamiltonian is
$$h_b = -\sum_{a,n} : T^a_n T^a_{-n} : = \frac{k+2}{2} d_b.$$
Now the lowest eigenvalue of $d_b$ is equal to the eigenvalue of the Casimir
operator $ - \sum_a T^a_0 T^a_0$ which is equal to $\frac12 j_0(j_0+1)$
where $j_0 =0,\frac12, 1, \dots k/2$ labels the vacuum represention
of $SU(2).$ Thus we obtain
$$ h_b -\frac12 j(j+1) \geq 0$$
for any $SU2)$ representation $j$ contained in $H_b.$

Similarly, on the fermionic sector $H_f$ we have
$$ d_f \geq \frac12 \ell(\ell+1)$$
since the vacuum eigenvalue of $d_f=
\frac{k+2}{2} h_f $ is zero; here $h_f = \frac{k+2}{8} n:\psi^a_n \psi^a_{-n}
: $ and $\frac12 \ell(\ell+1)$ is the eigenvalue of the invariant $-\sum_a
{{K'}^a}_0 {{K'}^a}_0$ in a given irreducible representation.
This inequality follows from the anticommutation relations of the fermion
operators ${\psi^a}_n:$ In order to increase the value of $\ell$ from zero
(in the vacuum) to a given value $\ell$ one must apply the fermion
operators at least for energies $n=1,2,\dots, \ell$ which leads  to the eigenvalue
$\frac12 \ell(\ell+1)$for $d_f.$

Thus we have
$$h_f = \frac{k+2}{2} d_f \geq \frac{k+2}{4} \ell(\ell+1).$$ 
Now by Lemma 3, in a given $(j,\ell)$ subrepresentation of the commuting
algebras $({T^a}_0)$ and $({{K'}^a}_0),$
$$\align {Q_+}^2 &= h_b + h_f -\frac12 (j+\ell)(j+\ell +1) \geq
h_f -\frac12 \ell(\ell+1) -j\ell \\ 
& \geq \frac{k+2}{4} \ell(\ell+1)-\frac12\ell(\ell+1)
-j\ell \geq  \ell( \frac{k}{4} (\ell +1) -j).\tag2.21 \endalign$$
In the subspace $H_{\Lambda}$ we have $\frac{(k+2)^2}{8} \geq h \geq h_b
\geq \frac12 j(j+1)$ so that $j\leq (k+1)/2.$ For $\ell\geq 2$ the right-hand-side
of (2.21) is strictly greater than zero. In the case $\ell=0$ we have $Q_+^2=
h_f +h_b - \frac12 j(j+1) -\frac18$ and the claim follows from the fact
that $h_b -\frac12 j(j+1)$ vanishes only on the vacuum sector. The remaining
case $\ell=1$ is clear from (2.21) if $j< k/2.$ But since we restrict to the
subspace $H_{\Lambda}$ where $h\leq (k+2)^2/8$ the cases $j \geq k/2$ are
excluded by the energy inequalities
$$h= h_b + h_f +\frac18 \geq \frac12 j(j+1) + \frac{k+2}{4} \ell(\ell+1)+
\frac18 = \frac12 j(j+1) +\frac{k+2}{2} +\frac18$$
for $\ell=1.$ \enddemo

\proclaim{Theorem 2} The family of operators $Q_A$ defined by the weight 
$(k,j_0)$ of a highest weight representation of $LG$ defines an
element in $K^1_{G}(G, k+2)$  for $G=SU(2).$ The value of the
invariant $I$ mod $k+2$  for this K theory class is equal to $2j_0+1$ and
therefore they are inequivalent for the allowed values $2j_0=
0,1,2,\dots, k.$ 
\endproclaim 

\demo{Proof} By (2.18) the value of the invariant $I(g)$ is given as 
$$I(g)= -\frac{i}{2} \tr_{\ker\, Q_+} \Gamma\tag2.22$$ 
since $Q_+$ anticommutes with $\Gamma.$ But since the kernel of $Q_+$
is equal to $H_0$ and $\Gamma=\psi^1_0\psi^2_0\psi^3_0 =
\sigma_1\sigma_2\sigma_3= i$ on the vacuum sector, we get $I(g)=
\frac12 \text{dim}\, H_0= 2j_0+1,$ where we have taken into account
that the dimension of $H_{f,vac}$ is two.  In particular, it follows 
that the trivial one dimensional representation $j_0=0$ gives 
the generator in $\Bbb Z/(k+2) \Bbb Z.$

\enddemo
\bf Remark \rm The construction of the operators $Q_A$ works for any
semisimple compact group $G,$ [M]. However, the twisted K theory
classes are not parametrized by a single invariant $I(g).$ Instead, 
one should study reductions of the K theory classes to various 
$SU(2)$ subgroups corresponding to a choice of simple roots of $G.$

\vskip 0.3in
\bf References \rm 

\vskip 0.3in

[AtSe] M.F. Atiyah and G. Segal: Twisted K theory, in  preparation.  

[AS] M.F. Atiyah and I. Singer: Index theory for skew-adjoint Fredholm 
operators.  I.H.E.S. Publ. Math. \bf 37, \rm 305 (1969)  

[BCMMS] P. Bouwknegt, A.L. Carey, V. Mathai, M. Murray: Twisted
K-theory  and K-theory of bundle gerbes. hep-th/0106194.

[FHT1] D. Freed, M. Hopkins, and C. Teleman: Twisted equivariant
K-theory with complex coefficients. math.AT/0206257. 

[FHT2] D. Freed, M. Hopkins, and C. Teleman: Twisted K theory and loop group
representations. math.AT/0312155.   
   
[K] V. Kac: \it Infinite Dimensional Lie Algebras. \rm Third
edition. Cambridge University Press, Cambridge (1990)

[KT] V. Kac and I. Todorov: Superconformal current algebra and their
unitary representations. Commun. Math. Phys. \bf 102, \rm 337 (1985)  

[M] J. Mickelsson:  Gerbes, (twisted) K-theory, and the supersymmetric
WZW model. hep-th/0206139. To be publ. in the proceedings of ``La
70eme Rencontre entre Physiciens Theoriciens et Mathematiciens'' in 
Strasbourg, May 23-25, 2002. 

[P] R. Palais:  On the homotopy type of certain groups of operators. Topology 3, 271--279
(1965)
  
\enddocument